\documentclass[amscd,amssymb,11pt]{amsart}

\usepackage[all]{xy}

\newtheorem{thm}{Theorem}[section]
\newtheorem{lem}[thm]{Lemma}
\newtheorem{cor}[thm]{Corollary}
\newtheorem{prop}[thm]{Proposition}

\theoremstyle{definition}

\newtheorem{ex}[thm]{Example}

\theoremstyle{remark}

\newtheorem{rem}[thm]{Remark}

\numberwithin{equation}{section}

\newcommand{\thmref}[1]{Theorem~\ref{#1}}
\newcommand{\corref}[1]{Corollary~\ref{#1}}
\newcommand{\secref}[1]{\S\ref{#1}}

\newcommand{\propref}[1]{Proposition~\ref{#1}}
\newcommand{\lemref}[1]{Lemma~\ref{#1}}

\newcommand{\hocolim}{\operatorname*{hocolim}}
\newcommand{\holim}{\operatorname*{holim}}

\newcommand{\totfib}{\operatorname*{TotFib}}
\newcommand{\totcofib}{\operatorname*{TotCofib}}

\newcommand{\Map}{\operatorname{Map}}
\newcommand{\MapS}{\operatorname{Map_{\mathcal S}}}

\newcommand{\A}{{\mathcal  A}}
\newcommand{\B}{{\mathcal  B}}
\newcommand{\C}{{\mathcal  C}}
\newcommand{\Sp}{{\mathcal  S}}

\newcommand{\PP}{{\mathcal P}}
\newcommand{\LL}{{\mathcal L}}

\newcommand{\Z}{{\mathbb  Z}}

\newcommand{\R}{{\mathbb  R}}

\newcommand{\Sinfty}{\Sigma^{\infty}}
\newcommand{\Oinfty}{\Omega^{\infty}}

\newcommand{\sm}{\wedge}

\newcommand{\ra}{\rightarrow}
\newcommand{\xra}{\xrightarrow}

\newcommand{\xla}{\xleftarrow}

\begin{document}

\title[Tate cohomology and periodic localization]{Tate cohomology and periodic localization of polynomial functors}             

\author[Kuhn]{Nicholas J.~Kuhn}                                  
\address{Department of Mathematics \\ University of Virginia \\ Charlottesville, VA 22903}    
\email{njk4x@virginia.edu}
\thanks{This research was partially supported by a grant from the National Science Foundation}     
 
\date{July, 2003.}


\begin{abstract}   
In this paper, we show that Goodwillie calculus, as applied to functors from stable homotopy to itself, interacts in striking ways with chromatic aspects of the stable category.

Localized at a fixed prime $p$, let $T(n)$ be the telescope of a $v_n$ self map of a finite $S$--module of type $n$.  The Periodicity Theorem of Hopkins and Smith implies that the Bousfield localization functor associated to $T(n)_*$ is independent of choices.  

Goodwillie's general theory says that to any homotopy functor $F$ from $S$--modules to $S$--modules, there is an associated tower under $F$, $\{P_dF\}$, such that $F \ra P_dF$ is the universal arrow to a $d$--excisive functor.

Our first main theorem says that $P_dF \ra P_{d-1}F$ always admits a homotopy section after localization with respect to $T(n)_*$ (and so also after localization with respect to  Morava $K$--theory $K(n)_*$).  Thus, after periodic localization, polynomial functors split as the product of their homogeneous factors.
  
This theorem follows from our second main theorem which is equivalent to the following: for any finite group $G$, the Tate spectrum $t_G(T(n))$ is weakly contractible.  This strengthens and extends previous theorems of Greenlees--Sadofsky, Hovey--Sadofsky, and Mahowald--Shick.  The Periodicity Theorem is used in an essential way in our proof.

The connection between the two theorems is via a reformulation of a result of McCarthy on dual calculus.
\end{abstract}
                                                  
\maketitle

\section{Introduction and main results} \label{introduction}

Over the past twenty years, beginning with the Nilpotence and Periodicity Theorems of E. Devanitz, M. Hopkins, and J. Smith \cite{dhs,hs,ravenel}, there has been a steady deepening of our understanding of stable homotopy as organized by the chromatic, or periodic, point of view.  During this same period, there have been many new results in homotopical algebra, many following the conceptual model offered by T. Goodwillies's calculus of functors \cite{goodwillie1, goodwillie2, goodwillie3}.

Here, and in a previous paper \cite{k2}, I prove theorems illustrating a beautiful interaction between these two strands of homotopy theory.  These results say that certain homotopy functors, stratified via Goodwillie calculus, decompose into their homogeneous strata, after periodic localization.  The first paper concerned a highly stuctured splitting of the important functor $\Sinfty \Oinfty$.  Ignoring the extra structure, one is left with an illustration of the main result here: after Bousfield localization with respect to a periodic homology theory, {\em all} polynomial endofunctors of stable homotopy split into a product of their homogeneous components.

We now explain our main results in more detail. 

The periodic homology theories we consider are $K(n)_*$, the $n^{th}$ Morava $K$--theory at a fixed prime $p$ and with $n > 0$, and the `telescopic' variants $T(n)_*$, where $T(n)$ denotes the telescope of a $v_n$--self map of a finite complex of type $n$.  A consequence of the Periodicity Theorem is that the associated Bousfield class $\langle T(n) \rangle$ is independent of the choice of both the complex and self map.  Also, we recall that $T(n)_*$--acyclics are $K(n)_*$--acyclic\footnote{The Telescope Conjecture asserts the converse, and, these days, is considered unlikely to hold for $n \geq 2$.}; thus the associated localization functors are related by $L_{K(n)} \simeq L_{K(n)}L_{T(n)}$.

Our use of concepts from Goodwillie calculus and localization theory require that we work within a good model category with homotopy category equivalent to the standard stable homotopy category.  Thus we work within the category $\Sp$, the category of $S$--modules of \cite{ekmm}.   

Goodwillie's general theory  then says that a homotopy functor $F: \Sp \ra \Sp$
 admits a universal tower of fibrations under $F$,
\begin{equation*}
\xymatrix{
&&& \vdots \ar[d] \\
&&& P_2F(X) \ar[d]^{p_2} \\
&&& P_1F(X) \ar[d]^{p_1} \\
F(X) \ar[rrr]^{e_0}  \ar[urrr]^{e_1} \ar[uurrr]^{e_2} &&& P_0F(X), 
}
\end{equation*}
such that \\

\noindent (1) $P_dF$ is $d$--excisive, and \\

\noindent (2) $e_d: F \ra P_dF$ is the universal natural transformation to a $d$-excisive functor.  \\

Our splitting theorem then is as follows. \\

\begin{thm} \label{splitting theorem}  Let $F: \Sp \ra \Sp$ be any homotopy functor.  For all primes $p$, $n \geq 1$, and $d \geq 1$, the natural map 
$$ p_d(X): P_dF(X) \ra P_{d-1}F(X)$$
admits a natural homotopy section after applying $L_{T(n)}$. \\
\end{thm}

The theorem can be reformulated as follows.  Let $D_dF(X)$ be the fiber of $p_d(X): P_dF(X) \ra P_{d-1}F(X)$.  Then $D_dF$ is both $d$--excisive and homogeneous: $P_{d-1}D_dF \simeq *$.  The theorem is equivalent to the statement that there is a natural weak equivalence of filtered spectra
$$ L_{T(n)}P_dF(X) \simeq \prod_{c=0}^d L_{T(n)}D_cF(X).$$

\begin{ex}  Here is the simplest example illustrating our theorem.  Let $p=2$.  For $k \in \Z$, let $\R P^{\infty}_k$ be the Thom spectrum of $k$ copies of the canonical line bundle over $\R P^{\infty}$.  \cite[Ex.5.7]{k1} implies that the cofibration sequence
\begin{equation} \label{cofib eq} \R P_{-1}^{\infty} \ra \R P_0^{\infty} \ra S^0
\end{equation}
splits after $K(n)$--localization, for all $n$, even though the connecting map $\delta: S^0 \ra \Sigma \R P_{-1}^{\infty}$ is nonzero in mod 2 homology.

As was, in essence, observed in a 1983 paper by J.Jones and S.Wegmann \cite{jw}, (\ref{cofib eq}) is the suspension of the special case $X = S^{-1}$ of a natural cofibration sequence of functors
\begin{equation} \label{cofib eq 2} (X \sm X)_{h\Z/2} \ra P_2(X) \ra X.
\end{equation}
One can also construct this sequence using Goodwillie calculus: see \secref{Z/p section}.

\thmref{splitting theorem} says that (\ref{cofib eq 2}) splits after applying $L_{T(n)}$ for all $n$ and $X$, even though the connecting map
$$ \delta: X \ra \Sigma (X \sm X)_{h\Z/2}$$
is often nontrivial before localization. \\
\end{ex}

\begin{rem}  There are various sorts of polynomial functors studied in the literature differing slightly from Goodwillie's $d$--excisive functors: R.McCarthy has studied $d$--additive functors \cite{mccarthy}, and his student A.Mauer--Oats \cite{mauer-oats} has studied an infinite family interpolating between additive and excisive.  As will be explained more fully in \secref{splitting thm section}, the analogue of \thmref{splitting theorem} holds in all these generalized settings. \\
\end{rem}

\begin{rem}  \thmref{splitting theorem} and \corref{Tate cor} below also has consequences for using the tower $\{P_dF(X)\}$ to understand $E_n^*(F(X))$, where $E_n$ is the usual p-complete integral height $n$ complex oriented commutative $S$--algebra.  Since it is known \cite{hovey} that $K(n)_*(X) = 0$ if and only if $E_n^*(X)=0$, our theorem says that the spectral sequence associated to the tower will collapse at $E_1$. \\
\end{rem}

\thmref{splitting theorem} is deduced from a rather different result in equivariant stable homotopy theory that we now describe.  

If $G$ is a finite group, let $G$--$\Sp$ denote the category of $S$--modules with $G$--action: the category of so--called `naive $G$--spectra'.  Note that any $S$--module can be considered as an object in $G$--$\Sp$ by giving it trivial $G$--action. 

For $Y \in G$--$\Sp$, we let $Y_{hG}$ and $Y^{hG}$ respectively denote associated homotopy orbit and homotopy fixed point $S$--modules.  There are various constructions in the literature, more \cite{gm} or less \cite{acd, ahearnkuhn, klein1,ww1} sophisticated, of a natural `Norm' map
$$ N(Y): Y_{hG} \ra Y^{hG}$$
satisfying the key property that $N(Y)$ is an equivalence if $Y$ is a finite free $G$--CW spectrum.  Let the Tate spectrum $t_G(Y)$ be defined as the cofiber of $N(Y)$.  As recently observed by J.Klein \cite{klein2}, up to weak equivalence, these constructions are unique: see \secref{Tate section}.

We prove the following vanishing theorem. \\

\begin{thm} \label{Tate theorem}  For all finite groups $G$, primes $p$, and $n \geq 1$, 
$$L_{T(n)}t_G (L_{T(n)}S) \simeq *.$$ 
\end{thm}

This theorem will turn out to be equivalent to the following corollary. \\

\begin{cor} \label{Tate corollary} If $T(n)$ is the telescope of any $v_n$--self map of a type $n$ complex, then $t_G (T(n)) \simeq *$. \\
\end{cor}

Besides implying \thmref{splitting theorem}, \thmref{Tate theorem} also leads to the following splitting result. \\

\begin{cor} \label{Tate cor}  For any $Y \in G$--$\Sp$, the fundamental cofibration sequence
$$ Y_{hG} \xra{N(Y)} Y^{hG} \ra t_G(Y)$$
splits after applying $L_{T(n)}$ for any $n$. \\
\end{cor}

One also immediately deduces results similar to \cite[Cor.8.7]{hoveystrickland}. \\

\begin{cor} For all finite groups $G$, the norm map induces an isomorphism
$$  T(n)_*(BG) \xra{\sim} T(n)^{-*}(BG).$$
Similarly, $L_{T(n)}(\Sinfty BG_+)$ is self dual in the category of $T(n)$--local spectra. \\
\end{cor}

Our two theorems are supported by three propositions. 

The first of these is a slight variant of results of R. McCarthy in \cite{mccarthy}, and establishes the connection between our two theorems.

We need to recall Goodwillie's classification of homogeneous polynomial functors \cite{goodwillie3}.  Let $\Sigma_d$ denote the $d^{th}$ symmetric group.  If our original functor $F$ is {\em finitary} (terminology from \cite{goodwillie3}), i.e. ~commutes with directed homotopy colimits, then $D_dF(X)$ is weakly equivalent to a homotopy orbit spectrum of the form  $$(C_F(d) \sm X^{\sm d})_{h \Sigma_d},$$
where $C_F(d)\in \Sigma_d$--$\Sp$ is determined naturally by $F$.  Important to us is that, even without the finitary hypothesis, there is a natural weak equivalence of the form
$$D_dF(X) \simeq (\Delta_dF(X))_{h \Sigma_d},$$
where $\Delta_dF$ is a functor determined naturally by $F$, taking values in the category  $\Sigma_d$--$\Sp$. \\

\begin{prop} \label{pullback prop}  Let $F: \Sp \ra \Sp$ be any homotopy functor.  For all $d \geq 1$, there is a homotopy pullback diagram
\begin{equation*}
\xymatrix{
P_dF(X) \ar[d]^{p_d} \ar[r] &
(\Delta_dF(X))^{h\Sigma_d} \ar[d]  \\
P_{d-1}F(X) \ar[r] &
t_{\Sigma_d}(\Delta_dF(X)).
}
\end{equation*}
This diagram is natural in both $X$ and $F$. \\
\end{prop}

Our other two propositions together imply \thmref{Tate theorem}. The first is a new very general observation about Tate spectra. \\

\begin{prop} \label{Tate prop}  Let $R$ be a ring spectrum and $E_*$ a homology theory.  If $t_{\Z/p}(R)$ is $E_*$--acyclic for all primes $p$, then so is $t_G(M)$ for all $R$--modules $M$ and for all finite groups $G$. \\
\end{prop}

We remark that, by standard arguments, $t_{\Z/p}(R) \simeq *$, and thus is certainly $E_*$--acyclic, for all primes $p$ such that $R_*$ is uniquely $p$--divisible.  In particular, to apply the proposition to the pair $(R,E_*) = (L_{T(n)}S, T(n)_*)$, one need to only look at the single prime involved in the periodic theory.

It is in proving our last proposition that deep results in periodic stable homotopy will be used. \\

\begin{prop} \label{Z/p prop} For all primes $p$ and $n \geq 1$, $L_{T(n)}t_{\Z/p}(L_{T(n)}S) \simeq *$. \\
\end{prop}

At this point we need to comment on results like \thmref{Tate theorem} in the literature.  

The main theorem of the 1988 article by M.Mahowald and P.Shick \cite{ms} can be restated as
\begin{equation} \label{ms equation} t_{\Z/2}(T(n)) \simeq *.
\end{equation}
A proof along their lines can presumably be done at odd primes as well.  We will see that the generalization of their theorem to all primes is equivalent to \propref{Z/p prop},  yielding one possible proof of that result.  We will offer a rather different proof, using the telescopic functors of Bousfield and the author \cite{bousfield1, k1, bousfield2}.

The main theorem of the 1996 article by J.Greenlees and H.Sadofsky \cite{gs} reads
\begin{equation} \label{gs equation} t_G(K(n)) \simeq *.
\end{equation}
Their proof is elementary (in the sense that consequences of the Nilpotence Theorem are not needed), but heavily uses two special facts about $K(n)$: it is complex oriented, and $K(n)_*(B\Z/p)$ is a finitely generated $K(n)_*$--module.  Note that neither of these two facts is available when considering $T(n)_*$.  For readers interested in the simplest proof of (\ref{gs equation}), it is hard to imagine improving upon the clever argument given in \cite[Lemma 2.1]{gs} showing that $t_{\Z/p}(K(n)) \simeq *$ , but our \propref{Tate prop} offers an alternative way to proceed starting from this.

The most substantial part of the main theorem of \cite{hsadofsky} says that  \begin{equation} \label{hs equation} L_{K(n)}t_G(L_{K(n)}S) \simeq *.
\end{equation}  Note that, were the Telescope Conjecture true, then (\ref{hs equation}) and \thmref{Tate theorem} would be equivalent; at any rate, the latter implies the former.  The authors prove their theorem by starting from (\ref{gs equation}), and then using the Periodicity Theorem, together with the technical heart of Hopkins and D. ~Ravenel's proof \cite{ravenel} that $L_{E(n)}$ is a smashing localization. Our proof of \thmref{Tate theorem} bypasses the need for the Hopkins--Ravenel argument.  \\

The rest of the paper is organized as follows.  In \secref{Tate section}, we review properties of the norm map and $t_G$, leading to a proof of \propref{Tate prop}.  In \secref{Z/p section}, supported by the appendix, we first discuss models for $L_Et_{\Z/p}(L_ES)$ for a general spectrum $E$, and then use telescopic functors to show that the model is contractible when $E = T(n)$.  The results of the previous two sections are combined in \secref{Tate thm section} yielding a proof of \thmref{Tate theorem}. Also in this section is a discussion of the equivalence of \thmref{Tate theorem} and \corref{Tate corollary}, with arguments similar in spirit to ones in  \cite{ms,hsadofsky}.  In \secref{polynomial section}, we review what we need to about $d$--excisive functors, and prove \propref{pullback prop} with arguments similar to those in \cite{mccarthy}.  In \secref{splitting thm section}, we prove our splitting results, \thmref{splitting theorem} and \corref{Tate cor}.

As is already evident, if $E$ is an $S$--module, we let $L_E$ denote Bousfield localization with respect to the associated homology theory $E_*$.  Throughout we also use the following conventions regarding functors taking values in $\Sp$.  We write $F \xra[\sim]{f} G$ if $f(X): F(X) \ra G(X)$ is a weak equivalence for all $X$.  By a weak natural transformation $f:F \ra G$ we mean a pair of natural tranformations of the form $F \xla[\sim]{g} H \xra{h} G$ or $F \xra{h} H \xla[\sim]{g} G$.  Finally, we say that a diagram of weak natural transformations commutes if, after evaluation on any object $X$, the associated diagram commutes in the stable homotopy category. \\

\noindent {\bf Acknowledgements}  I would like to thank various people who have helped me with this project.  Randy McCarthy and Greg Arone have helped me learn about Goodwillie towers.  Obviously Randy's paper \cite{mccarthy} has been important to my thinking, and Greg suggested the compelling reformulation of Randy's results given in \propref{pullback prop}.  Neil Strickland alerted me to the fact a conjecture of mine, that (\ref{hs equation}) was true, was already a theorem in the literature, and Hal Sadofsky similarly told me about Mahowald and Shick's theorem (\ref{ms equation}).  Our main results have been reported on in various seminars and conferences, e.g. at the A.M.S. meetings in January, 2003, and in Oberwolfach in March, 2003.

\section{Tate spectra and \propref{Tate prop}} 
\label{Tate section}

\subsection{Homotopy orbit and fixed point spectra}  For $G$ a fixed finite group, and $Y \in G$--$\Sp$, the $S$--modules $Y_{hG}$ and $Y^{hG}$ are defined in the usual way:
$$ Y_{hG} = (EG_+ \sm Y)/G, \text{\hspace{.1in} and \hspace{.1in}} Y^{hG} = (\Map_{\Sp}(EG_+,Y))^G.$$

Both of these functors take weak equivalences and cofibration sequences in $G$--$\Sp$ to weak equivalences and cofibration sequences in $\Sp$. (See \cite[Part I]{gm} for these sorts of facts.)

$Y_{hG}$ has an important additional property not shared with $Y^{hG}$: it commutes with filtered homotopy colimits.  

We record the following well known facts, which are fundamental when one considers the behavior of  $Y_{hG}$ and $Y^{hG}$ under Bousfield localization. \\

\begin{lem} \label{orbit lemma} If $f:Y \ra Z$ is a map in $G$--$\Sp$ that is an $E_*$--isomorphism, then $f_{hG}: Y_{hG} \ra Z_{hG}$ is also an $E_*$--isomorphism. \\
\end{lem}

\begin{lem} If $Y \in G$--$\Sp$ is $E$--local, so is $Y^{hG}$.
\end{lem}

\subsection{A characterization of the norm map}  A recent paper by Klein \cite{klein2} exploring axioms for generalized Farrell--Tate cohomology  leads to a nice characterization of norm maps, and thus Tate spectra. \\

\begin{prop} \label{norm char prop}Let $N_G(Y),N_G^{\prime}(Y): Y_{hG} \ra Y^{hG}$ be natural transformations such that both $N_G(\Sinfty G_+)$ and $N_G^{\prime}(\Sinfty G_+)$ are weak equivalences.  Then there is a unique weak natural equivalence $f(Y): Y_{hG} \xra{\sim} Y_{hG}$ such that the diagram
\begin{equation*}
\xymatrix{
Y_{hG} \ar[dr]_{f(Y)} \ar[r]^{N_G(Y)} &  Y^{hG}   \\
& Y_{hG} \ar[u]_{N_G^{\prime}(Y)} 
}
\end{equation*}
commutes.  It follows that the cofibers of $N_G(Y)$ and $N_G^{\prime}(Y)$ are naturally weakly equivalent.
\end{prop}

We sketch the proof, using the sorts of arguments in \cite{klein2}.  

Call a homotopy functor $H: \text{$G$--$\Sp$} \ra \Sp$ {\em homological} if it preserves homotopy pushout squares and filtered homotopy colimits.  Then Klein, in the spirit of \cite{ww2}, observes that any homotopy functor $F: \text{$G$--$\Sp$} \ra \Sp$ admits a universal left approximation by a homological functor, i.e. there exists homological functor $F^{hom}$, and a natural transformation $ F^{hom}(Y) \ra F(Y)$ satisfying the expected universal property.

Applying this to the case $F(Y) = Y^{hG}$, and observing that $Y_{hG}$ is homological, shows that there is a unique weak natural transformation $g:Y_{hG} \ra Y^{hG,hom}$ yielding a commutative diagram of weak natural transformations 
\begin{equation*}
\xymatrix{
Y_{hG} \ar[dr]_{g(Y)} \ar[r]^{N_G(Y)} &  Y^{hG}   \\
& Y^{hG,hom}. \ar[u] 
}
\end{equation*}
The right upward map is certainly an equivalence for $Y = \Sinfty G_+$, and, by assumption, so is the top map.  Thus $g$ is a weak natural transformation between homological functors that is an equivalence when $Y=\Sinfty G_+$. It follows that $g$ is weak equivalence.

Applying this same argument to $N_G^{\prime}$ yields the proposition.

\subsection{Tate spectra} We refer to any natural transformation $N_G$ as in the last proposition as a norm map.  The cofiber of $N_G(Y)$ is the associated Tate spectrum, denoted $t_G(Y)$.  Both $N_G$ and $t_G$ are unique in the sense of \propref{norm char prop}; their existence is shown in the various references cited in the introduction.

It is immediate that $t_G$ preserves weak equivalences and cofibration sequences.  

From \cite[Prop. I.3.5]{gm}, we deduce \\

\begin{lem} If $R$ is a (homotopy) ring spectrum with trivial $G$ action, and $M$ is an $R$--module, then $t_G(R)$ is a ring spectrum, and $t_G(M)$ is a $t_G(R)$--module.  Furthermore, $R^{hG} \ra t_G(R)$ is a map of $R$--algebra spectra. \\
\end{lem}

Fix $Y \in G$--$\Sp$.  For each subgroup $H$ of $G$, $Y$ can be regarded as being in $H$--$\Sp$ by restriction.  From \cite[pp.28--29]{gm}, one deduces \\

\begin{lem} The assignment $ G/H \longmapsto t_H(Y)$
defines a Mackey functor to the stable homotopy category.  Furthermore,
$ Y^{hH} \ra t_H(Y)$
is a map of Mackey functors. \\
\end{lem}

In \secref{polynomial section}, we will use the following familiar property of the norm map.  In the literature, this explicitly appears, with a short axiomatic proof, as \cite[Prop.2.10]{ahearnkuhn}. \\

\begin{lem} \label{map lemma} If $K$ is a finite free $G$--CW complex, then for all $Y \in G$--$\Sp$, $t_G(\MapS(K,Y)) \simeq *$.  \\
\end{lem}

\subsection{Proof of \propref{Tate prop}}

Recall that $R$ is a ring spectrum, and we are assuming that $t_{\Z/p}(R)$ is $E_*$--acyclic. We wish to show that $t_G(M)$ is also $E_*$ acyclic, for all $R$--modules $M$, and for all $G$.

We first note that we can assume $M = R$.  For $t_G(M)$ is a $t_G(R)$--module, and thus the former will be $E_*$--acyclic if the latter is.

Next we show that we can reduce to the case when $G$ is a $p$--group.  For each prime $p$ dividing the order of $G$, let $G_p < G$ be a $p$--Sylow subgroup.  Then we have \\

\begin{lem}  Given $Y \in G$--$\Sp$ and $E_*$ a generalized homology theory, $t_G(Y)$ will be $E_*$--acyclic if $t_{G_p}(Y)$ is $E_*$--acyclic for all $p$ dividing the order of $G$.
\end{lem}

\begin{proof} We recall that the completion of the Burnside ring $A(H)$ is denoted $\widehat{A}(H)$.  The assignment $G/H \longmapsto Y^{hH}$ is then an $\widehat{A}$--module Mackey functor in the sense of \cite{mm}. Thus so is $G/H \longmapsto t_H(Y)$, and then also $G/H \longmapsto E_*(t_H(Y))$.  Now \cite[Cor.4]{mm} implies the lemma. \\
\end{proof}

Having reduced \propref{Tate prop} to the case when $G$ is a $p$--group, and is thus solvable, the next lemma implies the proposition.

\begin{lem}  Let $K$ be a normal subgroup of $G$, $Q = G/K$, $R$ a ring spectrum, and $E_*$ a homology theory.  If $t_K(R)$ and $t_Q(R)$ are both $E_*$--acyclic, so is $t_G(R)$.
\end{lem}

\begin{proof}  For $Y \in G$--$\Sp$, consider the composite
$$ Y_{hG} \simeq (Y_{hK})_{hQ} \xra{N_K(Y)_{hQ}} (Y^{hK})_{hQ} \xra{N_Q(Y^{hK})} (Y^{hK})^{hQ} \simeq Y^{hG}.$$
We will know that this composite can be considered a norm map if we check that each of these maps is an equivalence when $Y = \Sinfty G_+$.  

As there is an equivalence of $S$--modules with $K$--action
$$\Sinfty G_+ \simeq \bigvee_{gK \in Q} \Sinfty K_+,$$
it follows that $N_K(\Sinfty G_+)$, and thus $N_K(\Sinfty G_+)_{hQ}$, is an equivalence.  

As there are equivalences of $S$--modules with $Q$--action
$$(\Sinfty G_+)^{hK} \xla[\sim]{N_K(\Sinfty G_+)} (\Sinfty G_+)_{hK} \xra[\sim]{} \Sinfty Q_+,$$
it follows that $N_Q((\Sinfty G_+)^{hK})$ is an equivalence.

We conclude from this discussion that if both $N_K(R)_{hQ}$ and $N_Q(R^{hK})$ are $E_*$--isomorphisms, then $N_G(R)$ will also be an $E_*$--isomorphism, and thus $t_G(R)$ will be $E_*$--acyclic. 

By assumption, $t_K(R)$ is $E_*$--acyclic.  Thus $N_K(R)$ is an $E_*$--isomorphism.  By \lemref{orbit lemma}, $N_K(R)_{hQ}$ is also.

By assumption, $t_Q(R)$ is $E_*$--acyclic.  As $t_Q(R^{hK})$ is a $t_Q(R)$--module, we conclude that $t_Q(R^{hK})$ is also $E_*$--acyclic, so that $N_Q(R^{hK})$ is an $E_*$--isomorphism.
\end{proof}

\section{Telescopic functors and  \propref{Z/p prop}} 
\label{Z/p section}

The goal of this section is to prove that $L_{T(n)}t_{\Z/p}L_{T(n)}S \simeq *$. 
We will prove this by establishing that the localized unit map
$$L_{T(n)}S \ra L_{T(n)}t_{\Z/p}L_{T(n)}S$$
is null.  

In outline our argument showing this is as follows.  It is well known that $t_{\Z/p}S$ can be written as certain inverse limit of Thom spectra.  Starting from this, we will show that the unit map $S \ra t_{\Z/p}S$ factors though an inverse limit of `connecting maps' associated to the Goodwillie tower of the functor $\Sinfty \Oinfty$ applied to spheres in negative dimensions.  We warn the reader of technical complications: odd primes are less pleasant than $p=2$, we use the theorems of W.H.Lin and J.Gunawardena establishing the Segal conjecture for $\Z/p$, and a key homological calculation is deferred to an appendix.

It will follow that the localized unit will factor through the inverse limit of the localized connecting maps.  That this inverse limit is null will then be an easy consequence of constructions of Bousfield and the author \cite{bousfield1, k1, bousfield2} showing that $L_{T(n)}$ factors through $\Oinfty$.  These `telescopic' constructions heavily use the Periodicity Theorem of Hopkins and Smith \cite{hs}, and thus are also heavily dependent on the Nilpotence Theorem of \cite{dhs}.

\subsection{Models for $L_Et_{\Z/p}L_EY$ and $L_Et_{\Sigma_p}L_EY$}

If $\alpha$ is an orthogonal real representation of a finite group $G$, we let $S(\alpha)$ and $S^{\alpha}$ respectively denote the associated unit sphere and one point compactified sphere.  Thus $S(\alpha)$ has an unbased $G$--action while the $G$--action on $S^{\alpha}$ is based, and there is a cofibration sequence of based $G$--spaces
$$ S(\alpha)_+ \ra S^0 \ra S^{\alpha}.$$

Fix a prime $p$, and let $\rho$ denote $\Sigma_p$ acting on $\R^p/\Delta(\R)$ in the usual way.  The action of $\Z/p < \Sigma_p$ on $S(\rho)$ is free, and one concludes that the infinite join $S(\infty \rho)$ is a model for $E\Z/p$.  This quickly leads to the following well known description of $t_{\Z/p}$.

\begin{lem} (Compare with \cite[Thm.16.1]{gm}.)  For $Y \in G$--$\Sp$, there is a natural weak equivalence
$$ t_{\Z/p}Y \simeq \holim_k \Sigma \MapS(S^{k\rho}, Y)_{h\Z/p}.$$
\end{lem}

We need a generalization of this.   \\

\begin{lem}  For $Y \in G$--$\Sp$, there is a natural weak equivalence
$$ L_Et_{\Z/p}L_EY \simeq \holim_k \Sigma L_E(\MapS(S^{k\rho}, Y)_{h\Z/p}).$$
If $(p-1)!$ acts invertibly on $E_*$, e.g. if $E$ is $p$--local, there is a natural weak equivalence
$$ L_Et_{\Sigma_p}L_EY \simeq \holim_k \Sigma L_E(\MapS(S^{k\rho}, Y)_{h\Sigma_p}).$$
These equivalences are also natural with respect to the partially ordered set of Bousfield classes $\langle E \rangle$, and there are commutative diagrams
\begin{equation*}
\xymatrix{
L_Et_{\Sigma_p}L_EY \ar[d] \ar[r]^-{\sim} &
\holim_k \Sigma L_E(\MapS(S^{k\rho}, Y)_{h\Sigma_p}) \ar[d]  \\
L_Et_{\Z/p}L_EY \ar[r]^-{\sim} &
\holim_k \Sigma L_E(\MapS(S^{k\rho}, Y)_{h\Z/p}).
}
\end{equation*}
\end{lem}

\begin{proof}  By definition, $L_E t_{\Z/p} L_E Y$ is the cofiber of 
$$ L_E N_{\Z/p}(L_E Y): L_E (L_E Y)_{h \Z/p} \ra L_E (L_E Y)^{h \Z/p}.$$
The domain of this map can be simplified:
$$ L_E Y_{h\Z/p} \ra L_E (L_EY)_{h\Z/p}$$
is an equivalence.  Meanwhile, the range of this map rewritten via the following chain of natural weak equivalences:
\begin{equation*}
\begin{split}
L_E (L_EY)^{h\Z/p} & \xla{\sim} (L_EY)^{h\Z/p} \\
  & \xra{\sim} \MapS(S(\infty \rho)_+, L_EY)^{h\Z/p} \\
  & \xra{\sim} \holim_k \MapS(S(k \rho)_+, L_EY)^{h\Z/p} \\
  & \xra{\sim} \holim_k L_E \MapS(S(k \rho)_+, L_EY)^{h\Z/p} \\
  & \xla{\sim} \holim_k L_E \MapS(S(k \rho)_+, L_EY)_{h\Z/p} \\
  & \xla{\sim} \holim_k L_E \MapS(S(k \rho)_+, Y)_{h\Z/p}. \\
\end{split} 
\end{equation*}
The crucial second to last equivalence here is induced by norm maps which are equivalences since $\Z/p$ acts freely on $S(k\rho)$.

Thus $L_E t_{\Z/p} L_E Y$ has been identified:
\begin{equation*}
\begin{split}
L_E t_{\Z/p} L_E Y & 
\xra{\sim} \holim_k \text{ cofiber } \{ L_E Y_{h\Z/p} \ra L_E 
\MapS(S(k\rho)_+, Y)_{h\Z/p} \} \\
  & \xra{\sim} \holim_k \Sigma L_E \MapS(S^{k\rho}, Y)_{h\Z/p}.
\end{split} 
\end{equation*}

The proof of the statements for $t_{\Sigma_p}$ are similar, noting that, under the hypothesis that $(p-1)!$ acts invertibly on $E_*$, the norm maps
$$ \MapS(S(k \rho)_+, L_EY)_{h\Sigma_p} \ra \MapS(S(k \rho)_+, L_EY)^{h\Sigma_p}$$
will still be equivalences.
\end{proof}

For $r \geq 0$, and $X$ an $S$--module, we let $D_rX = (X^{\sm r})_{h \Sigma_r}$, and we recall that there are natural transformations $\Sigma D_r X \ra D_r \Sigma X$.  Specializing to $r=p$, a quick check of definitions verifies the next lemma. \\

\begin{lem} There is a natural weak equivalence
$$ \Sigma^k D_p \Sigma^{-k}X \simeq \MapS(S^{k\rho}, X^{\sm p})_{h\Sigma_p},$$
and thus there is a $p$--local  equivalence
$$ t_{\Sigma_p}S \simeq \holim_k \Sigma^{k+1} D_pS^{-k}.$$
\end{lem}

Define $d_k: S \ra \Sigma^{k+1}D_p S^{-k}$ to be be the composite
\begin{equation*} \label{dk def}S \xra{\text{unit}} t_{\Sigma_p} S \longrightarrow \Sigma^{k+1}D_p S^{-k}.
\end{equation*}
As the restriction map $t_{\Sigma_p}S \ra t_{\Z/p}S$ is unital, our various observations combine to yield the following proposition. \\

\begin{prop} \label{vanishing prop} If $E$ is $p$--local, $L_E t_{\Z/p} L_ES \simeq *$ if and only if 
$$\holim_k L_E d_k: L_E S \ra \holim_k L_E \Sigma^{k+1} D_pS^{-k}$$
is null.
\end{prop}

\subsection{The Goodwillie tower of $\Sinfty \Oinfty$}

Recall that $\Sinfty Z$ denotes the suspension spectrum of a space $Z$, and that  $\Sinfty$ has right adjoint $\Oinfty$, where $\Oinfty X$ is the $\text{zero}^{th}$ space of a spectrum $X$.  

Let $P_r(X)$ denote the $r^{th}$ functor in the Goodwillie tower of the functor
$$\Sinfty \Oinfty: \Sp \ra \Sp.$$
Thus this Goodwillie tower has the form
\begin{equation*}
\xymatrix{
&&& \vdots \ar[d] \\
&&& P_3(X) \ar[d]^{p_3(X)} \\
&&& P_2(X) \ar[d]^{p_2(X)} \\
\Sinfty \Oinfty X \ar[rrr]^{e_1(X)}  \ar[urrr]^{e_2(X)} \ar[uurrr]^{e_3(X)} &&& P_1(X). 
}
\end{equation*}

This tower has the following fundamental properties. \\

\noindent{(1)} If $X$ is 0--connected, then $\Sinfty \Oinfty X \ra \holim_r P_r(X)$ is an equivalence. \\

\noindent{(2)} The fiber of $p_r(X): P_r(X) \ra P_{r-1}(X)$ is naturally weakly equivalent to $D_r(X)$. \\

\noindent{(3)} There are equivalences $D_1X \simeq P_1X \simeq X$, and via the second of these, $e_1(X): \Sinfty \Oinfty X \ra P_1X$ can be identified the with evaluation map $\epsilon(X): \Sinfty \Oinfty X \ra X$. \\

All of these properties can be deduced from Goodwillie's general theory.  For an explicit discussion of these (and more) see \cite{ahearnkuhn} or \cite{k2}.

\subsection{Telescopic functors}

Bousfield and the author have deduced the the following consequence of the Periodicity Theorem.

\begin{thm}  There exists a functor $\Phi_n: \text{Spaces} \ra \text{$S$--modules}$ and a natural weak equivalence $\Phi_n \Oinfty X \simeq L_{T(n)}X$.
\end{thm}

With the result stated at the level of homotopy categories, and with $K(n)$ replacing $T(n)$, this is the main theorem of \cite{k1}.  However the sorts of constructions given there, and in \cite{bousfield1} (for $n=1$), yield the theorem as stated: see \cite{bousfield2}.

This has the following immediate corollary \cite{k1,k2}. \\

\begin{cor} \label{Tn lifting cor} 
There is a natural factorization by weak $S$--module maps
\begin{equation*} 
\xymatrix{   & L_{T(n)}\Sinfty \Oinfty X \ar[dr]^{L_{T(n)}\epsilon(X)} &    \\
L_{T(n)}X \ar@{=}[rr] \ar[ur]^{\eta_n(X)}  &  &  L_{T(n)}X. \\
}
\end{equation*}
\end{cor}

To use this, we recall an observation about {\em reduced} homotopy functors,  functors $F:\Sp \ra \Sp$ such that $F(X)$ is contractible whenever $X$ is.  Goodwillie observes that then there is an induced weak natural transformation
$$ \Sigma F(X) \longrightarrow F(\Sigma X).$$
The naturality is with respect to both $X$ and $F$.
For example, if $F = D_r$, this natural transformation agrees with the one discussed previously.

In particular, we can apply this construction to both the domain and range of the natural transformation
$$ L_{T(n)}P_p \ra L_{T(n)}P_1,$$
evaluated on $\Sigma^{-k}X$ for all $k \geq 0$.
Recalling that $L_{T(n)}$ commutes with suspension and $P_1(X) \simeq X$, we obtain maps
$$ \holim_k \Sigma^k L_{T(n)}P_p(\Sigma^{-k}X) \ra L_{T(n)}X.$$

\begin{thm} \label{vanishing thm}  $\displaystyle \holim_k \Sigma^k L_{T(n)}P_p(\Sigma^{-k}X) \ra L_{T(n)}X$ admits a homotopy section.
\end{thm}

\begin{proof} A section is given by $\displaystyle \holim_k \Sigma^k (L_{T(n)}e_p(\Sigma^{-k}X) \circ \eta_n(\Sigma^{-k}X))$.
\end{proof}

\subsection{Specialization to odd spheres}

Standard homology calculations as in \cite{clm,bmms} imply the next lemma. \\

\begin{lem}  Localized at an odd prime $p$, $D_r S^k \simeq *$ for odd $k \in \Z$, and for $2 \leq r \leq p-1$.  Thus (for all primes $p$) the natural map
$$ \holim_k \Sigma^k P_{p-1}(S^{-k}) \ra S$$
is a $p$--local equivalence.
\end{lem}

Continuing the cofibration sequence $ D_p X \ra P_p(X) \ra P_{p-1}(X)$
one step to the right defines a natural transformation 
$$\delta(X): P_{p-1}(X) \ra \Sigma D_p X.$$
Localized at $p$, define $\delta_k: S \ra \Sigma^{k+1}D_p S^{-k}$ to be the composite
\begin{equation*} \label{delta k def}S \xla{\sim} \Sigma^kP_{p-1}(S^{-k})  \xra{\Sigma^k \delta(S^{-k})} \Sigma^{k+1}D_p S^{-k}.
\end{equation*}

\begin{prop} \label{null prop} $\displaystyle \holim_k L_{T(n)} \delta_k: L_{T(n)} S \ra \holim_k L_{T(n)} \Sigma^{k+1} D_pS^{-k}$
is null.
\end{prop}
\begin{proof}  Localized at $p$, there is a cofibration sequence
$$ \holim_k L_{T(n)} \Sigma^{k} P_pS^{-k} \ra L_{T(n)} S \ra \holim_k L_{T(n)} \Sigma^{k+1} D_pS^{-k}.$$
\thmref{vanishing thm} says that the first map has a section.  Thus the second map is null.
\end{proof}

\subsection{Proof of \propref{Z/p prop}}
A comparison of \propref{vanishing prop} with \propref{null prop} shows that we will have proved \propref{Z/p prop} once we check the following lemma. \\

\begin{lem} \label{d = delta lemma} $\displaystyle \holim_k  d_k: S \ra \holim_k \Sigma^{k+1} D_pS^{-k}$
factors through \\ 
$\displaystyle \holim_k \delta_k:  S \ra \holim_k  \Sigma^{k+1} D_pS^{-k}$.
\end{lem}

\begin{proof} W.H.Lin's theorem \cite{lin}, when $p=2$, and J.Gunawardena's theorem \cite{g, agm}, when $p$ is odd, can be stated in the following way: 
$$\displaystyle \holim_k  d_k: S \ra \holim_k \Sigma^{k+1} D_pS^{-k}$$
is $p$--adic completion.  It follows that we need to check that 
$\displaystyle \holim_k  \delta_k \in \pi_0(\holim_k \Sigma^{k+1} D_pS^{-k}) \simeq \Z_p$
is a topological generator.  As topological generators of $\Z_p$ are detected mod $p$, the next lemma, whose proof is deferred to the appendix, completes our argument.
\end{proof}

\begin{lem} \label{homology lem} $\delta(S^{-1}): P_{p-1}(S^{-1}) \ra  \Sigma D_p S^{-1}$ is nonzero in mod $p$ homology.
\end{lem}

\section{The proofs of \thmref{Tate theorem} and \corref{Tate corollary}}
\label{Tate thm section}

We begin this section by noting how \propref{Z/p prop} and \propref{Tate prop} together imply \thmref{Tate theorem}.  \propref{Z/p prop} can be restated as saying that $t_{\Z/p}L_{T(n)}S$ is $T(n)_*$--acyclic.  Recalling that the localization of a ring spectrum (e.g.~ $S$) is again a ring spectrum, \propref{Tate prop} can then be applied to the pair $(R,E_*) = (L_{T(n)}S, T(n)_*)$, to conclude that $t_{G}L_{T(n)}S$ is $T(n)_*$--acyclic for all $G$.  This is a restatement of \thmref{Tate theorem}.

Now we turn to showing how \corref{Tate corollary} can be deduced from \thmref{Tate theorem}, and vice versa.

We need to review some of the fine points of the Periodicity Theorem. (A good reference for this is \cite{ravenel}.) We fix a prime $p$, and work with $p$--local spectra. A finite spectrum $F$ is of type $n$ if $K(n)_*(F) \neq 0$, but $K(i)_*(F)=0$ for $i < n$.  Let $\C_n = \{ \text{finite }F \ | \ F \text{ has type at least } n \}$.  Then every $F \in \C_n$ admits a $v_n$ self map: a map $f: \Sigma^d F \ra F$ such that $K(n)_*(f)$ is an isomorphism, but $K(i)_*(f) = 0$ for all $i \neq n$.  If $n > 0$, then $d$ will necessarily be positive.  In all cases, $f$ is unique and natural up to iteration.  Thus there is a well defined functor from $\C_n$ to spectra sending $F$ to $v_n^{-1}F$, the telescope of any $v_n$ self map of $F$.  We note that $v_n$ preserves both cofibration sequences and retracts.

The Thick Subcategory Theorem says that any thick subcategory of the category of $p$--local spectra, i.e. ~any collection of $p$--local finite spectra closed under cofibration sequences and retracts, is $\C_n$ for some $n \geq 0$.

We recall that $L_{T(n)}$ denotes $L_{v_n^{-1}F}$ for any $F$ of type $n$.  From the facts stated above, it is easily verified that this is independent of choice of $F$, and that for all $F \in \C_n$, $L_{T(n)}(F) = v_n^{-1}F$. Finally we note that if $F$ has type $n$ and $F^{\prime}$ has type $i \neq n$, then $v_n^{-1}F \sm v_i^{-1}F^{\prime} \simeq *$.  \\

\begin{lem} Fix a finite group $G$.  The following conditions are equivalent. \\

\noindent (1) $t_G(L_{T(n)}S)$ is $T(n)_*$--acyclic. \\

\noindent (2) For all $F \in \C_n$, $t_G(v_n^{-1}F) \simeq *$. \\

\noindent (3) For all type $n$ complexes $F$, $t_G(v_n^{-1}F) \simeq *$. \\

\noindent (4) There exists a type $n$ complex $F$ such that $t_G(v_n^{-1}F) \simeq *$. \\
\end{lem}

Note that statement (1) is the conclusion of \thmref{Tate theorem} and (3) is the conclusion of \corref{Tate corollary}.

Clearly (2) implies (3), which in turn implies (4).  To see that (4) implies (2), note that the collection of $F \in \C_n$ such that $t_G(v_n^{-1}F) \simeq *$ forms 
a thick subcategory contained in $\C_n$.  Such a thick subcategory will equal all of $\C_n$ if it contains {\em any} type $n$ finite.  (This type of reasoning appears in \cite{ms}.)

Now suppose (1) holds.  Since $v_n^{-1}F \simeq L_{T(n)}F$, it is an $L_{T(n)}S$--module, and we see that $t_G(v_n^{-1}F)$ is $(v_n^{-1}F)_*$--acyclic for all finite $F$ of type $n$.  It is easy to find a type $n$ finite $F$ that is a ring spectrum; thus so is $R = v_n^{-1}F$.  But then $t_G(R)$ will be an $R_*$--acyclic $R$--module, and thus contractible, i.e. statement (4) holds.

It remains to show that (2) implies (1).  We reason as in \cite{hsadofsky}.  

Define finite spectra $F(0), \dots, F(n)$ by first setting $F(0) = S$, and then recursively defining $F(i+1)$ to be the cofiber of a $v_i$ self map of $F(i)$.  

Ravenel \cite{rav1} observes that if $f: \Sigma^d X \ra X$ is a self map with cofiber $C$ and telescope $T$, then $\langle X \rangle = \langle C \vee T \rangle$.  Applying this $n$ times leads to an equality of Bousfield classes
$$ \langle S \rangle = \langle F(n) \vee \bigvee_{i=0}^{n-1} v_i^{-1}F(i) \rangle.$$

Smashing this with $t_G(L_{T(n)}S)$, and noting that $$t_G(L_{T(n)}S) \sm F(n) \simeq t_G(L_{T(n)}F(n)) \simeq t_G(v_n^{-1}F(n)),$$ leads to 
$$ \langle t_G(L_{T(n)}S) \rangle = \langle t_G(v_n^{-1}F(n)) \vee \bigvee_{i=0}^{n-1} t_G(L_{T(n)}S) \sm v_i^{-1}F(i) \rangle.$$

Smashing this with $T(n)$, and noting that $T(n) \sm v_i^{-1}F(i) \simeq *$ if $i < n$, leads to 
$$ \langle T(n) \sm t_G(L_{T(n)}S) \rangle = \langle T(n) \sm t_G(v_n^{-1}F(n)) \rangle.$$
If (2) holds, then the right side of this last equation is the Bousfield class of a contactible spectrum.  Thus so is the left, i.e. ~(1) holds.

\section{Polynomial functors and Tate cohomology}
\label{polynomial section}

In this section we sketch a proof of \propref{pullback prop}.  As I hope will be clear, this proposition is just a variant of \cite[Prop.4]{mccarthy}, and our proof uses precisely the same ideas that McCarthy does.

\subsection{Review of Goodwillie calculus}
In the series of papers \cite{goodwillie1, goodwillie2, goodwillie3}, Tom Goodwillie has developed his theory of polynomial resolutions of homotopy functors.  We need to summarize some aspects of Goodwillie's work as they apply to functors from $S$--modules to $S$--modules.  Throughout we cite the version of \cite{goodwillie3} of June, 2002.

In \cite{goodwillie2}, Goodwillie begins by defining and studying the {\em total homotopy fiber} of a cubical diagram.  For example the total homotopy fiber of a square
\begin{equation*}
\xymatrix{
X_0 \ar[d] \ar[r] &
X_1 \ar[d]  \\
X_2 \ar[r] &
X_{12}
}
\end{equation*}
is the homotopy fiber of the evident map from $X_0$ to the homotopy pullback of the square with $X_0$ omitted.  A cubical diagram is then {\em homotopy cartesian} if its total fiber is weakly contractible.  Dual constructions similarly define {\em total homotopy cofibers} and {\em homotopy cocartesion cubes}.  We note that in a stable model category like $\Sp$, a cubical diagram is homotopy cartesian exactly when it is homotopy cocartesion.

A cubical diagram is {\em strongly homotopy cocartesion} if each of its 2 dimensional faces is homotopy cocartesion.  A functor is then said to be {\em $d$--excisive} if it takes strongly homotopy cocartesion $(d+1)$--cubical diagrams to homotopy cartesian cubical diagrams.  

In \cite{goodwillie3}, given a functor $F$, Goodwillie proves the existence of a tower $\{P_dF\}$ under $F$ so that $F \ra P_dF$ is the universal arrow to a $d$--excisive functor, up to weak equivalence.   

For functors with range in a stable model category, Goodwillie \cite{goodwillie3} gives a description of how $D_dF(X)$, the fiber of $P_dF(X) \ra P_{d-1}F(X)$, can be computed by means of cross effects.  We describe how this goes in our setting. 

Let $ F: \Sp \ra \Sp$  be a functor.  Let ${\mathbf  d} = \{1,2,\dots, d\}$.  In \cite[\S 3]{goodwillie3}, $cr_dF$, the $d^{th}$ {\em cross effect} of $F$, is defined to the the functor of $d$ variables given as the total homotopy fiber
\begin{equation*}  
(cr_dF)(X_1, \dots, X_d) = \totfib_{T \subset \mathbf d} F(\bigvee_{i \in {\mathbf d} - T} X_i).
\end{equation*}

A $d$--variable homotopy functor $H: \Sp^d \ra \Sp$ is {\em reduced} if ~$H(X_1, \dots, X_d)$ is contractible whenever any of the $X_i$ are.  Given such a functor, its {\em multilinearization} $\LL(H): \Sp^d \ra \Sp$ is defined by the formula
\begin{equation} \label{linearization}
 \LL(H)(X_1, \dots, X_d) = \hocolim_{n_i \ra \infty} \Omega^{n_1 + \dots + n_d} H(\Sigma^{n_1}X_1, \dots, \Sigma^{n_d}X_d).
\end{equation}
This will be 1-excisive in each variable.

Now define $\Delta_dF: \Sp \ra \Sigma_d$--$\Sp$ by the formula 
$$ \Delta_dF(X) = \LL(cr_dF)(X, \dots, X).$$

Then \cite[Theorems 3.5, 6.1]{goodwillie3} says that there is a natural weak equivalence
\begin{equation} \label{Dd formula}
 D_dF(X) \simeq (\Delta_dF)(X)_{h\Sigma_d}.
\end{equation}

We need to explain some of the ideas behind this formula.  

Firstly, $\Delta_d(F) \ra \Delta_d(P_dF)$ is always an equivalence, and it follows that one can assume the original functor $F$ is $d$--excisive.  

If $F$ is $d$--excisive then $cr_dF$ is already 1--excisive in each variable \cite[Prop.3.3]{goodwillie3}, and so $\Delta_dF(X)$ can be identified with $(cr_dF)(X, \dots, X)$. In this case, the natural map 
$$ D_dF(X) \ra P_dF(X)$$
identifies with the natural transformation
$$ \alpha_d(X): (\Delta_dF)(X)_{h\Sigma_d} \ra F(X)$$
defined to be the composite
$$ (\Delta_dF)(X)_{h\Sigma_d} \ra F(\bigvee_{i=1}^d X)_{h\Sigma_d} \ra F(X).$$
Here the second map is induced by the fold map $ \bigvee_{i=1}^d X \ra X$.

Goodwillie proves (\ref{Dd formula}) by verifying that $cr_d(\alpha_d)$ is an equivalence, so that $D_d(\alpha_d)$ is an equivalence.  Enroute to this, he shows that there is a natural equivariant weak equivalence 
$$ cr_d(\Delta_d F) \simeq \Sigma_{d+} \sm cr_d F.$$

\subsection{Dual constructions}

In \cite{mccarthy}, McCarthy investigates `dual calculus'.  In this spirit,  replacing wedges by products, fibers by cofibers, etc., leads to constructions dual to the above.  In particular, given $F: \Sp \ra \Sp$, we define $cr^dF: \Sp^d \ra \Sp$ by the formula
\begin{equation*}  
(cr^dF)(X_1, \dots, X_d) = \totcofib_{T \subset \mathbf d} F(\prod_{i \in T} X_i),
\end{equation*}
and then we define $\Delta^dF: \Sp \ra \Sigma_d$--$\Sp$ by
$$ \Delta^dF(X) = \LL(cr^dF)(X, \dots, X).$$

Because both the domain and range of $F$ is a stable model category, one sees that each of the natural transormations
$$cr_dF \ra cr^dF$$
and 
$$ \Delta_dF \ra \Delta^dF$$
are weak equivalences.

If $F$ is $d$--excisive then $\Delta^dF(X)$ can be identified with $(cr^dF)(X, \dots, X)$.  In this case, we define the weak natural transformation
$$ \alpha^d(X): F(X) \ra (\Delta_dF)(X)^{h\Sigma_d}$$
to be the zig--zag composite
$$F(X)  \ra F(X^d)^{h\Sigma_d} \ra (\Delta^dF)(X)^{h\Sigma_d} \xla{\sim} (\Delta_dF)(X)^{h\Sigma_d}.$$
Here the first map is induced by the diagonal $ X \ra X^d$.

Arguments dual to Goodwillie's show that the next lemma holds. \\

\begin{lem} \label{dual calc lemma} (Compare with \cite[Lemmas 3.7,3.8]{mccarthy}.)  Let $F: \Sp \ra \Sp$ be $d$--excisive. \\

\noindent (1) $cr^d(\alpha^d)$, and thus $D_d(\alpha^d)$, is an equivalence. \\

\noindent (2) There is a natural equivariant weak equivalence 
$$ cr^d(\Delta^d F) \simeq \MapS(\Sigma_+, cr^d F).$$
\end{lem}
 
\subsection{Proof of \propref{pullback prop}}

\propref{pullback prop} is a formal consequence of \lemref{dual calc lemma}.  First of all, we observe the following.

\begin{lem} (Compare with \cite[proof of Prop.4]{mccarthy}.)  Let $F$ be $d$--excisive.  Then $t_{\Sigma_d}(\Delta_dF)$ is $(d-1)$--excisive.  Thus the cofibration sequence
$$ D_d((\Delta_dF)^{h\Sigma_d}) \ra P_d(\Delta_dF)^{h\Sigma_d} \ra P_{d-1}((\Delta_dF)^{h\Sigma_d})$$
identifies with the norm sequence
$$ (\Delta_dF)_{h\Sigma_d} \ra (\Delta_dF)^{h\Sigma_d} \ra t_{\Sigma_d}(\Delta_dF).$$
\end{lem} 
\begin{proof}  For the first statement, we check that $cr^d(t_{\Sigma_d}(\Delta_dF)) \simeq *$:
$$ cr^d(t_{\Sigma_d}(\Delta_dF)) \simeq t_{\Sigma_d}(cr^d(\Delta_dF)) \simeq t_{\Sigma_d}(\MapS(\Sigma_+, cr^d F)) \simeq *.$$
Here we have used \lemref{dual calc lemma}(2) and \lemref{map lemma}.

As $(\Delta_dF)_{h\Sigma_d}$ is $d$--excisive and homogeneous, the second statement follows.
\end{proof}

Now we turn to the proof of \propref{pullback prop}.  We can assume that $F$ is $d$--excisive.  Assuming this, the last lemma implies that the weak natural tranformation $\alpha^d(X): F(X) \ra (\Delta_dF(X))^{h\Sigma_d}$ induces a commutative diagram of weak natural transformations
\begin{equation*}
\xymatrix{
D_dF(X) \ar[d] \ar[r] & (\Delta_dF(X))_{h\Sigma_d} \ar[d]  \\
P_dF(X) \ar[d] \ar[r] & (\Delta_dF(X))^{h\Sigma_d} \ar[d]  \\
P_{d-1}F(X) \ar[r] &  t_{\Sigma_d}(\Delta_dF(X)).
}
\end{equation*}
In this diagram each of the vertical columns is a homotopy fibration sequence of $S$--modules.  The top map is a weak  equivalence thanks to \lemref{dual calc lemma}(1).  Thus the bottom square is a homotopy pullback diagram.

\subsection{Polynomial functor variants} \label{variants subsection} McCarthy and his student Mauer--Oats \cite{mauer-oats} have explored various different notions of what it might mean to say a  functor $F: \A \ra \B$ is polynomial of degree at most $d$, with $d$--excisive and $d$--additive as two special cases. In these variants $\B$ should surely be a reasonable model category, but $\A$ can often be a category with much less structure.  As a hint of why this might be true, note that the definition of cross effects only uses the existence of finite coproducts in $\A$.  

If $\B$ is any stable model category admitting norm maps, and $\A$ is also appropriately stable, then the evident analogue of \propref{pullback prop} still holds.  The discussion above goes through with one little change:  the formula (\ref{linearization}) for the (multi)linearization process $\LL$ needs to be adjusted to reflect the notion of degree 1 functor at hand.  Note that our proof of \propref{pullback prop} didn't use this formula (nor did McCarthy's arguments in \cite{mccarthy}).

Of relevance to the next section, we note that these variants of $\LL$ are still homotopy colimits, and thus preserve $E_*$--isomorphisms.

\section{Localization and the proofs of \thmref{splitting theorem} and \corref{Tate cor}}  \label{splitting thm section}

In this section, we show how our vanishing Tate cohomology result, \thmref{Tate theorem}, leads to the splitting results \thmref{splitting theorem} and \corref{Tate cor}.  To simplify notation, we let $L = L_{T(n)}$.

\begin{proof}[Proof of \corref{Tate cor}]
Let $Y$ be an $S$--module with $G$ action.  We wish to show that the norm sequence
$$ Y_{hG} \xra{N(Y)} Y^{hG} \ra t_G(Y)$$
splits after applying $L$. Thus we need to construct a left homotopy inverse to $L(N(Y))$. 

The localization map $Y \ra LY$ induces a commutative diagram 
\begin{equation*}
\xymatrix{
Y_{hG} \ar[d] \ar[rr]^-{N(Y)} && Y^{hG} \ar[d]   \\
(LY)_{hG}  \ar[rr]^-{N(LY)} && (LY)^{hG}.  \\
}
\end{equation*}
Applying $L$ to this yields the diagram
\begin{equation*}
\xymatrix{
L(Y_{hG}) \ar[d]^{\wr} \ar[rr]^-{L(N(Y))} && L(Y^{hG}) \ar[d]^{L(\eta^{hG})}   \\
L((LY)_{hG})  \ar[rr]^-{L(N(LY))}_-{\sim} && L((LY)^{hG}).  \\
}
\end{equation*}
Here the left vertical map is an equivalence, as homology isomorphisms are preserved by taking homotopy orbits (\lemref{orbit lemma}).  The lower map, $L(N(LY))$, is an equivalence by \thmref{Tate theorem}: its cofiber, $L(t_G(LY))$, is a module over $L(t_G(LS))$, and is thus contractible.  

Our desired left homotopy inverse is now obtained by composing the right vertical map of the diagram with the inverses of the two indicated equivalences.
\end{proof}

\begin{proof}[Proof of \thmref{splitting theorem}]  We are given a functor $F: \Sp \ra \Sp$ and wish to prove that 
$$ D_dF(X) \ra P_dF(X) \ra P_{d-1}F(X)$$
splits after applying $L$.  Thus we need to construct a left homotopy inverse to $LD_dF(X) \ra LP_dF(X)$. 

We need a lemma that plays the role that \lemref{orbit lemma} played in the previous proof.  Call a natural transformation $F \ra G$ an $E_*$--isomorphism, if $F(X) \ra G(X)$ is an $E_*$--isomorphism for all $X$.  

Formula (\ref{Dd formula}) says that $D_d$ is the composition of constructions each of which preserve $E_*$--isomorphisms, and thus we have 

\begin{lem} If $F \ra G$ is an $E_*$--isomorphism, then so is $D_dF \ra D_dG$. \\
\end{lem}
 
\begin{rem}  This lemma holds for the variants on the notion of $d$--excisive, as discussed above in \secref{variants subsection}.  \\
\end{rem}

Armed with this lemma, \thmref{splitting theorem} is proved as follows.

The localization natural transformation $F \ra LF$, together with \propref{pullback prop}, induce a commutative diagram 
\begin{equation*}
\xymatrix{
D_dF \ar[d] \ar[r] & D_d(LF) \ar[d] \ar[r]^-{\sim} & \Delta_d(LF)_{h\Sigma_d} \ar[d]  \\
P_dF \ar[r] & P_d(LF) \ar[r] &  \Delta_d(LF)^{h\Sigma_d}.
}
\end{equation*}
Applying $L$ to this, gives the diagram
\begin{equation*}
\xymatrix{
LD_dF \ar[d] \ar[r]^-{\sim} & LD_d(LF) \ar[d] \ar[r]^-{\sim} & L(\Delta_d(LF)_{h\Sigma_d}) \ar[d]^{\wr}  \\
LP_dF \ar[r] & LP_d(LF) \ar[r] &  L(\Delta_d(LF)^{h\Sigma_d}).
}
\end{equation*}
Here the top left natural transformation is an equivalence by the lemma just stated.  The right vertical natural transformation is an equivalence by \thmref{Tate theorem}, as its cofiber, $L(t_{\Sigma_d}(\Delta_d(LF))$, is an  $L(t_G(LS))$--module, when evaluated on any $X$. (Though not necessarily local, due to the hocolimit construction $\LL$, $\Delta_d(LF)(X)$ is nevertheless an $LS$--module.)

Our desired left homotopy inverse is now obtained by composing the natural transformation along the bottom of this diagram with the inverses of the three indicated equivalences.
\end{proof}

\appendix

\section{Proof of \lemref{homology lem}}

We begin with some needed notation.

Recall that $P_r$ denotes the $r^{th}$ Goodwillie approximation to the functor $\Sinfty \Oinfty$.  We let
$$\delta(X): P_{r-1}(X) \ra \Sigma D_r X$$ 
denote the connecting map for the cofibration sequence
$$ D_r X \ra P_r(X) \ra P_{r-1}(X).$$

Given any reduced homotopy functor $F: \Sp \ra \Sp$, we let 
$$ \Delta(X): \Sigma F(X) \ra F(\Sigma X)$$
denote the canonical natural map.

Fixing a prime $p$, all homology will be with $\Z/p$ coefficients.  The Steenrod operations act on $H_*(X)$ as operations lowering dimensions.  To unify the `even' prime and odd prime cases, we let $\PP^1 = Sq^2$, when $p=2$.  Thus, for all primes $p$, $\PP^1$ lowers degree by $2p-2$.

The goal of this appendix is to prove \lemref{homology lem}, which we restate more precisely.

\begin{lem} \label{homology lem 2} $\delta_*: H_{-1}(P_{p-1}(S^{-1})) \ra H_{-1}(\Sigma D_p(S^{-1}))$ is an isomorphism of one dimensional $\Z/p$--modules. \\
\end{lem}

Recall that $H_*(D_rX)$ is a known functor of $H_*(X)$, both additively, and as a module over the Steenrod algebra.  Furthermore, the behavior of $\Delta_*: H_*(\Sigma D_r X) \ra H_*(D_r \Sigma X)$ is known.  See \cite{clm, bmms}.

Naturality implies that there is a commutative diagram:
\begin{equation*}
\xymatrix{
H_{-1}(P_{p-1}(S^{-1})) \ar[d]_{\wr} \ar[rr]^-{\delta_*} &&  H_{-1}(\Sigma D_p(S^{-1})) \ar[d]_{\wr}  \\
H_{-1}(P_{p-1}(\Sigma^{-1}H\Z)) \ar[rr]^-{\delta_*} && H_{-1}(\Sigma D_p(\Sigma^{-1}H\Z)) \\
H_{2p-3}(P_{p-1}(\Sigma^{-1}H\Z))  \ar[u]_{\PP^1_*} \ar[d]^{\Delta_*} \ar[rr]^-{\delta_*} &&  H_{2p-3}(\Sigma D_p(\Sigma^{-1}H\Z)) \ar[u]^{\wr}_{\PP^1_*} \ar[d]_{\wr}^{\Delta_*} \\
H_{2p-3}(\Sigma^{-2} P_{p-1}(\Sigma H\Z)) \ar[rr]^-{\Sigma^{-2}\delta_*} && H_{2p-3}(\Sigma^{-1} D_p(\Sigma H\Z)),
}
\end{equation*}
where the top vertical maps are induced by the inclusion $S^{-1} \ra \Sigma^{-1}H\Z$.  The top square is a square of homology groups of lowest degree.  That the indicated maps are isomorphisms, all between one dimensional vector spaces, is an easy consequence of facts from \cite{clm, bmms}.  For example, the middle right map is an isomorphism due to the Nishida relation
$$ \PP^1_* \beta Q^1 x = \beta Q^0 x \in H_{-2}(D_p(\Sigma^{-1}H\Z)),$$
for $x \in H_{-1}(\Sigma^{-1}H\Z)$.

Using this diagram, to show that the top map is nonzero, and thus an isomorphism, it suffices to show that the lower left map and the bottom map are each isomorphisms.  We state each of these as a separate lemma (one in dual form). \\

\begin{lem} $\Delta_*: H_{2p-3}(P_{p-1}(\Sigma^{-1}H\Z)) \ra H_{2p-3}(\Sigma^{-2} P_{p-1}(\Sigma H\Z))$ is an isomorphism of one dimensional $\Z/p$--modules.
\end{lem}
\begin{proof}  When $p=2$, $\Delta$ is an equivalence, and so $\Delta_*$ is an isomorphism.

When $p$ is odd, the situation is more complicated, and we proceed as follows.  We have a commutative diagram
\begin{equation*}
\xymatrix{
H_{2p-3}(P_{p-1}(\Sigma^{-1}H\Z)) \ar[d]^{\wr} \ar[rr]^{\Delta_*} &&  H_{2p-3}(\Sigma^{-2} P_{p-1}(\Sigma H\Z)) \ar[dd]^{\wr}  \\
H_{2p-3}(P_{2}(\Sigma^{-1}H\Z)) \ar[d] && \\
H_{2p-3}(\Sigma^{-1}H\Z) \ar[rr]^{\Delta_*}_{\sim} &&  H_{2p-3}(\Sigma^{-1}H\Z)
}
\end{equation*}
with indicated isomorphisms.  Thus, to show the top map is an isomorphism, we need to check that the lower left map is an isomorphism.  Equivalently, we need to check that
$$ \delta_*: H_{2p-3}(\Sigma^{-1} H \Z) \ra H_{2p-3}(\Sigma D_2 \Sigma^{-1} H \Z)$$
is zero.  The map $\delta: \Sigma^{-1} H \Z \ra \Sigma D_2 \Sigma^{-1} H \Z$ factors through
$$ \Delta: \Sigma^2 D_2 \Sigma^{-2} H \Z \ra \Sigma D_2 \Sigma^{-1} H \Z,$$
and this map is zero on $H_{2p-3}$: the range is one dimensional, spanned by the suspension of a $*$--decomposable of the form $x*y$, with $x \in H_{-1}(\Sigma^{-1}H\Z)$ and $y \in H_{2p-3}(\Sigma^{-1}H\Z)$.  But nonzero $*$--decomposables are never in the image of $\Delta_*: H_*(\Sigma D_2(X)) \ra H_*(D_2(\Sigma X))$.
\end{proof}

With our final lemma, we have reached the heart of the matter. \\

\begin{lem} $\delta^*: H^{2p-1}(\Sigma D_p(\Sigma H\Z)) \ra H^{2p-1}(P_{p-1}(\Sigma H\Z))$ is an isomorphism of one dimensional $\Z/p$--modules.
\end{lem}
\begin{proof}  Since $\Sigma H \Z$ is 0--connected, the Goodwillie tower $P_r(\Sigma H\Z)$ converges strongly to $\Sinfty \Oinfty (\Sigma H \Z) = \Sinfty S^1$.  Thus the associated 2nd quadrant spectral sequence converges strongly to $H^*(S^1)$.  For this to happen, $\PP^1(x)$ must be in the image of $\delta^*$, where $x \in H^1(P_{p-1}(\Sigma H\Z))$ is a nonzero element, for otherwise $\PP^1(x) \neq 0 \in H^{2p-1}(S^1)$.
    
Thus $\delta^*$ is nonzero, and is thus an isomorphism. 
\end{proof}

\begin{rem}  In work in progress, the author is studying the spectral sequence converging to $H^*(\Oinfty X)$ with $E_1^{-r,*+r} = H^{*}(D_r X)$.  The sort of argument just given generalizes to show that the first interesting differential is $d_{p-1}: H^{*-1}(D_p X) \ra H^*(X)$.  This differential is determined by $H^*(X)$ as a module over the Steenrod algebra, and has image imposing the unstable condition on $H^*(X)$.
\end{rem}

\end{document}